\pdfoutput=1
\RequirePackage{ifpdf}
\ifpdf % We are running pdfTeX in pdf mode
\documentclass[pdftex]{sigma}
\else
\documentclass{sigma}
\fi

\numberwithin{equation}{section}
\newtheorem{Theorem}{Theorem}[section]
\newtheorem{Lemma}[Theorem]{Lemma}
\newtheorem{Proposition}[Theorem]{Proposition}
{\theoremstyle{definition}
\newtheorem{Definition}[Theorem]{Definition}
\newtheorem{Remark}[Theorem]{Remark}
}

\begin{document}

\allowdisplaybreaks

\renewcommand{\thefootnote}{$\star$}

\renewcommand{\PaperNumber}{033}

\FirstPageHeading

\ShortArticleName{Hyperk\"ahler Manifolds of Curves in Twistor Spaces}

\ArticleName{Hyperk\"ahler Manifolds of Curves in Twistor Spaces\footnote{This paper is a~contribution to the Special
Issue on Progress in Twistor Theory. The full collection is available at
\href{http://www.emis.de/journals/SIGMA/twistors.html}{http://www.emis.de/journals/SIGMA/twistors.html}}}

\Author{Roger BIELAWSKI}
\AuthorNameForHeading{R.~Bielawski}

\Address{Institut f\"ur Differentialgeometrie, Universit\"at Hannover, \\
Welfengarten 1, D-30167 Hannover, Germany}

\Email{\href{mailto:bielawski@math.uni-hannover.de}{bielawski@math.uni-hannover.de}}

\URLaddress{\url{http://www.ifam.uni-hannover.de/~bielawski/}}

\ArticleDates{Received November 06, 2013, in f\/inal form March 19, 2014; Published online March 28, 2014}

\Abstract{We discuss hypercomplex and hyperk\"ahler structures obtained from higher deg\-ree curves in complex spaces
f\/ibring over ${\mathbb{P}}^1$.}


\Classification{53C26; 53C28; 32G10; 14H15}

\renewcommand{\thefootnote}{\arabic{footnote}}
\setcounter{footnote}{0}

\section{Introduction}

This article is concerned with constructions of hypercomplex and hyperk\"ahler structures from curves of arbitrary
degree, and with their properties.
It has been motivated by three sources.
First and foremost, the work of Nash~\cite{Nash}, who gave a~new twistor construction of hyperk\"ahler metrics on moduli
spaces of ${\rm SU}(2)$ magnetic monopoles.
Second, the so-called generalised Le\-gendre transform construction of hyperk\"ahler metrics, due to Lindstr\"om and
Ro\v{c}ek~\cite{LR}, which often leads to curves of higher genus.
Third, the well-known fact that the smooth locus of the Hilbert scheme of (local complete intersection) curves of degree
$d$ and genus $g$ in ${\mathbb{P}}^3$ has, if nonempty, dimension $4d$.
Because of author's ${\mathbb{H}}$-bias, the factor $4$ seems to him to call for some sort of quaternionic structure.

It turns out that the above three situations can be put in a~common framework.
Let $M$ be a~connected hypercomplex or a~hyperk\"ahler manifold.
The twistor space of $M$ is a~complex manifold $Z$ f\/ibring over ${\mathbb{P}}^1$, $\pi:Z\to {\mathbb{P}}^1$, and
equipped with an antiholomorphic involution $\sigma$, which covers the antipodal map on ${\mathbb{P}}^1$.
The manifold $M$ is recovered as a~connected component of the space of $\sigma$-equivariant sections
$s:{\mathbb{P}}^1\to Z$ with normal bundle $N\simeq {\mathcal{O}}(1)^n$.
We now ask: {\em what happens if we consider $\sigma$-invariant curves of higher degree in $Z$?} It turns out that we
still obtain a~hypercomplex manifold, as long as we require that the normal bundle $N$ of such a~curve $C$ satisf\/ies the
``stability condition'' $H^\ast(N\otimes_{{\mathcal{O}}_C}\pi^\ast{\mathcal{O}}_{{\mathbb{P}}^1}(-2))=0$.
This is the condition shown by Nash~\cite{Nash} to hold for spectral curves of monopoles, and used by him to describe
the hypercomplex structure of monopole moduli spaces.
His argument works in the general situation considered here.
Moreover, this new hypercomplex manifold is (pseudo)-hyperk\"ahler if $M$ was a $4$-dimensional hyperk\"ahler manifold.

The hypercomplex manifolds obtained this way have interesting properties.
For example, they are biholomorphic, with respect to any complex structure $I_\zeta$, $\zeta\in {\mathbb{P}}^1$, to
(unramif\/ied covering of) an open subset of the {\em smooth locus} of the Hilbert scheme~$Z_\zeta^{[d]}$ of~$d$ points in
the f\/ibre $Z_\zeta=\pi^{-1}(\zeta)$ (here $d$ is the degree of curves under consideration).

Hyperk\"ahler monopole moduli spaces arise in the above manner from the twistor space of $S^1\times {\mathbb{R}}^3$.
If we consider instead the twistor space ${\mathbb{P}}^3-{\mathbb{P}}^1$ of the f\/lat ${\mathbb{R}}^4$, we shall obtain
hyperk\"ahler structures on manifolds parameterising curves in ${\mathbb{P}}^3$ not intersecting a~f\/ixed line.
In the simplest case, that of twisted normal cubics, the resulting $12$-dimensional metric is still f\/lat, and the
question arises what happens for other admissible values of genus and degree.
Equally interesting is the question what happens for twistor spaces of compact hyperk\"ahler or hypercomplex manifolds.

The dif\/ferential geometry of hyperk\"ahler manifolds obtained from higher degree curves is richer than just
hyperk\"ahler geometry.
This has been already observed in~\cite{BS2} in the case of monopole moduli spaces.
A f\/irst step in understanding this geometry is a~description of natural objects on such a~manifold~$M$ directly in terms
of the complex manifold $Z$ containing the higher degree curves, without passing to the usual (higher-dimensional)
twistor space of the hyperk\"ahler structure.
Here we give such a~description for hyperholomorphic connections on vector bundles on~$M$.
The novelty is that we construct canonical connections via a~canonical splitting on the level of sections of vector
bundles, without a~corresponding splitting of vector bundles on~$Z$.

Finally, we make a~technical remark.
If we want to consider higher degree curves, then the Kodaira moduli spaces (of complex submanifolds in a~complex
manifold) are not enough: such curves will almost certainly degenerate, while the hyperk\"ahler metric will remain
smooth.
In fact, even for the usual twistor spaces, one is often enough led to consider a~singular space~$Z$, which already
contains all the information needed, and it is then simpler to work directly with~$Z$ rather than resolving its
singularities.
For these reasons we replace throughout the Kodaira moduli spaces with the Douady space~${\mathcal{D}}(Z)$.
We recall the necessary def\/initions and facts in the next section.

\section{Some background material
\label{back}
}

We gather here some necessary facts and def\/initions from complex analysis.
A good reference is~\cite{GR}.

We work in the category of complex spaces, i.e.\
${\mathbb C}$-ringed spaces locally modelled on ${\mathbb C}$-ringed subspaces of domains $U$ of ${\mathbb C}^n$ def\/ined
by f\/initely many holomorphic functions in $U$.
In particular, we allow nilpotents in the structure sheaf.

For a~complex space $(X,{\mathcal{O}}_X)$, the topological (Urysohn--Menger), analytic (Chevalley) and algebraic (Krull)
notions of local dimension $\dim_xX$ coincide.
A complex space is called {\em reduced} if its structure sheaf has no nilpotents.
It is {\em pure dimensional} (or equidimensional) if $\dim_x X$ is the same at all points of $X$.
For a~$0$-dimensional space $X$, its {\em length} is the sum of dimensions of stalks ${\mathcal{O}}_{X,x}$.

The {\em cotangent sheaf}\, $\Omega^1_X$ of a~complex space $X$ is def\/ined by glueing together the cotangent sheaves
${\mathcal{I}}_V/{\mathcal{I}}_V^2$ of local model spaces $(V, {\mathcal{O}}_{{\mathbb C}^n}/{\mathcal{I}}_V)$.
The {\em tangent sheaf}\, ${\mathcal{T}}_X$ is its dual $\operatorname{Hom}(\Omega^1_X,{\mathcal{O}}_X)$.
The {\em tangent space}\, at $x\in X$ is $T_xX={\mathcal{T}}_{X,x}/\mathfrak{m}_x {\mathcal{T}}_{X,x}$, where
$\mathfrak{m}_x$ is the maximal ideal of the stalk ${\mathcal{O}}_{X,x}$.

If $A$ is a~local ring with maximal ideal $\mathfrak{m}$ and $M$ is a~module over $A$, then a~sequence
$(f_1,\dots,f_r)\in \mathfrak{m}$ is said to be {\em regular for $M$} if $f_{i+1}$ does not divide $0$ in
$M/(f_1M+\dots+f_iM)$ for $i=0,\dots, r-1$.
The length of a~maximal regular sequence is called the {\em depth of $M$}.
A~complex space $X$ is said to be {\em Cohen--Macaulay} if $\operatorname{depth}{\mathcal{O}}_{X,x}=\dim_xX$ for all
$x\in X$.
A~reduced $1$-dimensional complex space is always Cohen--Macaulay.

A closed complex subspace $X$ of a~complex space $Z$ is called {\em regularly embedded} of codimension~$r$ if every of
its local def\/ining ideals $ {\mathcal{I}}_{X,x}$ has depth~$r$.
If $Z$ is a~complex manifold, and $X\subset Z$ is regularly embedded, then~$X$ is called a~{\em local complete
intersection}, abbreviated l.c.i.
This is an intrinsic condition and it implies that $X$ is Cohen--Macaulay.

For a~closed complex subspace $X$ of a~complex space $Z$ with ideal sheaf $ {\mathcal{I}}_X\subset {\mathcal{O}}_Z$, the
normal sheaf of $X$ in $Z$ is ${\mathcal{N}}_{X/Z}=\operatorname{Hom}_{{\mathcal{O}}_X}( {\mathcal{I}}_X/
{\mathcal{I}}_X^2,{\mathcal{O}}_X)$.
If $X\subset Z$ is regularly embedded, then both ${\mathcal{N}}_{X/Z}$ and $ {\mathcal{I}}_X/ {\mathcal{I}}_X^2$ are
locally free.
If, moreover, $Z$ is smooth and $X$ is reduced, then the conormal sheaf $ {\mathcal{I}}_X/ {\mathcal{I}}_X^2$ f\/its into
an exact sequence
\begin{gather*}
0\to {\mathcal{I}}_X/ {\mathcal{I}}_X^2\to \Omega^1_Z|_X\to \Omega^1_X\to 0,
\end{gather*}
and consequently we have an exact sequence
\begin{gather*}
0\to {\mathcal{T}}_X\to {\mathcal{T}}_Z|_X\to {\mathcal{N}}_{X/Z}\to \mathcal{E}xt^1\big(\Omega^1_X,{\mathcal{O}}_X\big).
\end{gather*}

A celebrated theorem of Douady~\cite{Dou} states that for any complex space $Z$ there exists a~complex space
${\mathcal{D}}(Z)$ parameterising all pure-dimensional compact complex subspaces of $Z$.
In addition, there exists a~{\em universal family} on ${\mathcal{D}}(Z)$, i.e.\
a~complex subspace $Y$ of ${\mathcal{D}}(Z)\times Z$, which def\/ines a~double f\/ibration
\begin{gather}
{\mathcal{D}}(Z)\stackrel{\nu}{\longleftarrow} Y \stackrel{\mu}{\longrightarrow} Z
\label{seqq0}
\end{gather}
with the following properties:
\begin{itemize}\itemsep=0pt
\item[(i)] $\nu$ is f\/lat and proper; \item[(ii)] if $S$ is a~complex space, $T\subset S\times Z$ a~complex subspace with
properties stated in (i), then there exists a~unique holomorphic map $f:S\to {\mathcal{D}}(Z)$ such that $T\simeq
S\times_{{\mathcal{D}}(Z)} Y$.
\end{itemize}
The f\/ibre of $Y$ over any $m\in {\mathcal{D}}(Z)$ is the complex subspace $X$ of $Z$ corresponding to $m$ and the
restriction of the normal sheaf ${\mathcal{N}}_Y$ of $Y$ in ${\mathcal{D}}(Z)\times Z$ to this f\/ibre is
${\mathcal{N}}_{X/Z}$.

At any $X\in {\mathcal{D}}(Z)$, there is a~canonical isomorphism $T_X{\mathcal{D}}(Z)\simeq H^0(X,{\mathcal{N}}_{X/Z})$.
The Douady space ${\mathcal{D}}(Z)$ is smooth at $X$ if $\operatorname{Ext}^1_{{\mathcal{O}}_Z}(
{\mathcal{I}}_X,{\mathcal{O}}_X)=0$.
For a~regularly embedded $X$, this is equivalent to $H^1(X,{\mathcal{N}}_{X/Z})=0$.
For a~projective $Z$, ${\mathcal{D}}(Z)$ is the Hilbert scheme parameterising compact subschemes of $Z$.

\section{Hypercomplex and hyperk\"ahler structures\\ from higher degree curves}\label{hgc}

For us, a~{\em curve} means a~compact, pure $1$-dimensional, Cohen--Macaulay complex space.
If the curve is reduced, then the Cohen--Macaulay assumption is redundant.

Let $Z$ be a~pure $(n+1)$-dimensional complex space and $\pi:Z\to {\mathbb{P}}^1$ a~holomorphic surjection.
For a~closed subspace $X$ of $Z$ we continue to write $\pi$ for its restriction to $X$.
Each f\/ibre $X_\zeta$, i.e.\
the complex subspace $(\pi^{-1}(\zeta), {\mathcal{O}}_X/\pi^\ast\mathfrak{m}_\zeta$), is a~closed subspace of~$X$.
If ${\mathcal{F}}$ is a~sheaf on $X$, we write ${\mathcal{F}}(k)$ for ${\mathcal{F}}\otimes \pi^\ast
{\mathcal{O}}_{{\mathbb{P}}^1}(k)$.
Furthermore, we write
\begin{gather*}
H_X=\pi^\ast H^0({\mathbb{P}}^1,{\mathcal{O}}(1)) \subset H^0(X,{\mathcal{O}}_X(1)).
\end{gather*}

We consider the set ${\mathcal{X}}$ of curves $C$ in $Z$ such that:
\begin{itemize}\itemsep=0pt
\item
$\pi:C\to {\mathbb{P}}^1$ is f\/inite-to-one;
\item
the Douady space ${\mathcal{D}}(Z)$ is smooth at $C$;
\item
the sheaf ${\mathcal{N}}_{C/Z}(-2)$ is acyclic, i.e.\
$h^i({\mathcal{N}}_{C/Z}(-2))=0$ for $i=0,1$.
\end{itemize}
These conditions are open, so ${\mathcal{X}}$ is an open submanifold of ${\mathcal{D}}(Z)$.

A useful characterisation of {\em pure dimension + Cohen--Macaulay + finite cover of ${\mathbb{P}}^1$} is provided by

\begin{Lemma}\label{CM}
Let $\pi:X\to Y$ be a~finite-to-one closed surjective holomorphic map from a~complex space $X$ to a~complex manifold
$Y$.
Then the following conditions are equivalent:
\begin{itemize}\itemsep=0pt
\item[$(i)$] $X$ is pure dimensional and Cohen--Macaulay;
\item[$(ii)$] $\pi$ is flat;
\item[$(iii)$] $\pi_\ast{\mathcal{O}}_X$ is locally free;
\item[$(iv)$] all fibres of $\pi$ have the same length.
\end{itemize}
\end{Lemma}

\begin{proof}
Follows from~\cite[Corollary 18.17]{Eis1} and the equality $\dim_xX=\dim_{\pi(x)}Y+\dim_{x} X_{\pi(x)}$ for f\/lat
maps~\cite[Proposition II.2.11]{GR}.
\end{proof}

Thus, for a~$C\in {\mathcal{X}}$, each f\/ibre $C_\zeta$ is a~$0$-dimensional complex space of constant length $d$, which
we call the {\em degree of the curve $C$}.

For any $C\in {\mathcal{X}}$ and $\zeta\in {\mathbb{P}}^1$, consider the map ${\mathcal{N}}_{C/Z}(-2)\stackrel{\cdot
s}{\longrightarrow} {\mathcal{N}}_{C/Z}(-1)$, where $s\in H_C$ with $s(\zeta)=0$.
This map must be injective, since a~nontrivial kernel sheaf would be supported on a~$0$-dimensional subspace, and so it
would have nontrivial sections, which would then map to nontrivial sections of ${\mathcal{N}}_{C/Z}(-2)$.
We have, therefore, the exact sequence
\begin{gather}
0\to {\mathcal{N}}_{C/Z}(-2)\stackrel{\cdot s}{\longrightarrow} {\mathcal{N}}_{C/Z}(-1)\longrightarrow
{\mathcal{N}}_{C/Z}(-1)|_{C_\zeta}\to 0.
\label{seq}
\end{gather}
Taking cohomology and considering a~generic $\zeta$ shows that
\[
h^0( {\mathcal{N}}_{C/Z}(-1))=dn \qquad \text{and}\qquad h^1(
{\mathcal{N}}_{C/Z}(-1))=0.
\]
Similarly, taking a~section $t$ of ${\mathcal{O}}_{{\mathbb{P}}^1}(2)$ with zeros at $\zeta$ and $\tilde\zeta$ gives the
exact sequence
\begin{gather}
0\to {\mathcal{N}}_{C/Z}(-2)\stackrel{\cdot t}{\longrightarrow} {\mathcal{N}}_{C/Z}\longrightarrow
{\mathcal{N}}_{C/Z}|_{C_\zeta}\oplus {\mathcal{N}}_{C/Z}|_{C_{\tilde\zeta}}\to 0,
\label{seqq}
\end{gather}
from which we derive
\begin{gather}
h^0( {\mathcal{N}}_{C/Z})=2dn,
\qquad
h^1( {\mathcal{N}}_{C/Z})=0.
\label{dim}
\end{gather}
The def\/initions and~\eqref{dim} imply
\begin{Proposition}
The subset ${\mathcal{X}}_d$ of ${\mathcal{X}}$, consisting of curves of degree $d$, is open in ${\mathcal{D}}(Z)$ and
smooth of dimension~$2dn$.
\end{Proposition}

For regularly embedded subspaces some of the assumptions on $C$ are automatically fulf\/illed

\begin{Lemma}\label{ooo}
Let $C$ be a~regularly embedded compact subspace of $Z$ such that $\pi:C\to {\mathbb{P}}^1$ is finite-to-one and
$H^\ast({\mathcal{N}}_{C/Z}(-2))=0$.
Then $C\in {\mathcal{X}}$.
\end{Lemma}

\begin{proof}
The normal sheaf of a~regularly embedded subspace is locally free.
Hence,~\eqref{seq} implies that each f\/ibre $C_\zeta$ has length $h^0( {\mathcal{N}}_{C/Z}(-1))/n$.
Owing to Lemma~\ref{CM}, $C$ is equidimensional and Cohen--Macaulay.
Moreover, for a~regularly embedded $C$,~\eqref{dim} implies that ${\mathcal{D}}(Z)$ is smooth at~$C$.
\end{proof}

Since the linear system $H_C$ is base-free,~\eqref{seqq} implies that the following sequence is also exact:
\begin{gather*}
0\to {\mathcal{N}}_{C/Z}(-2)\to {\mathcal{N}}_{C/Z}(-1)\otimes H_C\to {\mathcal{N}}_{C/Z}\to 0,
%\label{seq-1}
\end{gather*}
and, consequently, there is a~canonical isomorphism
\begin{gather}
H^0({\mathcal{N}}_{C/Z})\simeq H^0({\mathcal{N}}_{C/Z}(-1))\otimes H_C.
\label{tensor}
\end{gather}
We denote by $E$ a~vector bundle over ${\mathcal{X}}_d$, the f\/ibre of which at $C$ is $H^0({\mathcal{N}}_{C/Z}(-1))$.
In the notation of~\eqref{seqq0}, $E= \nu_\ast\bigl({\mathcal{N}}_Y\otimes \mu^\ast{\mathcal{O}}_Z(-1)\bigr)$.
We also write $H$ for the trivial bundle of rank~$2$, with f\/ibre $H_C=\pi^\ast H^0({\mathbb{P}}^1,{\mathcal{O}}(1))$.
The decomposition~\eqref{tensor} induces a~decomposition (of holomorphic vector bundles)
\begin{gather*}
T{\mathcal{X}}_d\simeq E\otimes H.
%\label{decompose}
\end{gather*}
For every $\zeta\in {\mathbb{P}}^1$, we have a~subbundle $Q_\zeta$ of $T{\mathcal{X}}_d$ of rank $dn$ def\/ined as
\begin{gather*}
Q_\zeta=E\otimes s,
\end{gather*}
where $s$ is a~section of ${\mathcal{O}}_{{\mathbb{P}}^1}(1)$ vanishing at $\zeta$.
In other words, the f\/ibre of $Q_\zeta$ at $C$ consists of sections of ${\mathcal{N}}_{C/Z}$ vanishing on $C_\zeta$.

\begin{Proposition}\label{integr}
The distribution $Q_\zeta$ is integrable.
\end{Proposition}

\begin{proof}
Let $X$ and $Y$ be vector f\/ields with values in $Q_\zeta$, and $t\mapsto\gamma_t(m)$, $t\mapsto \delta_t(m)$ their
integral curves beginning at $m$.
The bracket of $[X,Y]$ at any $p\in {\mathcal{X}}_d$ can be computed as $\lim\limits_{t\to 0}\dot\nu(t)/t^2$, where
\begin{gather*}
\nu(t)=\delta_{-t}(\gamma_{-t}(\delta_t(\gamma_t(m)))).
\end{gather*}
Consider the corresponding deformations of curves in $Z$.
Since $X$ and $Y$ take values in $Q_\zeta$, the deformation $\nu(t)$ f\/ixes the f\/ibre $C_\zeta$, and hence $\dot\nu(t)\in
Q_\zeta$.
Thus $[X,Y]_m\in {Q_\zeta}|_m$.
\end{proof}

Now suppose that $Z$ is equipped with an antiholomorphic involution $\sigma$ which covers the antipodal map on
${\mathbb{P}}^1$.
The submanifold $({\mathcal{X}}_d)^\sigma$ of $\sigma$-invariant curves in ${\mathcal{X}}_d$ is either empty or of real
dimension $2dn$.
In the latter case, $({\mathcal{X}}_d)^\sigma$ is canonically a~hypercomplex manifold.
For each complex structure $I_\zeta$, $\zeta\in {\mathbb{P}}^1$, the bundle of $(0,1)$-vectors is $Q_\zeta$.
We state the result as follows.
\begin{Theorem}\label{hc}
Let $Z$ be an equidimensional complex space equipped with a~holomorphic surjection $\pi:Z\to {\mathbb{P}}^1$ and an
antiholomorphic involution $\sigma$ covering the antipodal map on ${\mathbb{P}}^1$.
Then the subset ${\mathcal{M}}_d$ of the smooth locus of ${\mathcal{D}}(Z)$, consisting of $\sigma$-invariant curves $C$
of degree~$d$ such that $H^\ast({\mathcal{N}}_{C/Z}(-2))=0$, is, if nonempty, a~hypercomplex manifold of real dimension
$2d(\dim Z -1)$.
\end{Theorem}

\begin{Remark}
[a~generalisation] One can consider a~more general situation.
Let ${\mathcal{O}}(1)$ a~globally generated line bundle on a~complex space $Z$ and let $H\subset
H^0(Z,{\mathcal{O}}(1))$ be a~base-free linear system.
Let $X$ be a~regularly embedded compact subspace $X$ of $Z$ such that $h^1({\mathcal{N}}_{X/Z})=0$, $H|_X$ is base-free,
and the kernel of the natural map ${\mathcal{N}}_{X/Z}(-1)\otimes H\to {\mathcal{N}}_{X/Z}$ has no cohomology in
dimensions $0$, $1$ and $2$.
We obtain again a~canonical isomorphism
\begin{gather*}
H^0(X,{\mathcal{N}}_{X/Z})\simeq H^0({\mathcal{N}}_{X/Z}(-1))\otimes H,
\end{gather*}
so that the open subset ${\mathcal{X}}\subset
{\mathcal{D}}(Z)$, consisting of such $X$, is a~manifold, the tangent bundle of which decomposes as
$T{\mathcal{X}}\simeq E\otimes H$.
Such a~decomposition is known under various names, in particular as an {\em almost Grassmann structure} or a~{\em conic
structure}~\cite{Ma}.
Once again, for every $s\in H$, we have a~subbundle $Q_s$ of $T{\mathcal{X}}$ (of rank $\dim {\mathcal{X}}/\dim H$), and
the proof of Proposition~\ref{integr} can be repeated to show that $Q_s$ is integrable.
We can also consider compatible real structures.
In particular, we did not need to restrict ourselves to $1$-dimensional subspaces of $Z$ in order to obtain hypercomplex
manifolds.
\label{Grass}
\end{Remark}

Returning to the hypercomplex manifold ${\mathcal{M}}_d$, recall that, for a~$C\in {\mathcal{M}}_d$, each f\/ibre
$C_\zeta=C\cap \pi^{-1}(\zeta)$ is a~$0$-dimensional complex space of length $d$, and so, for each $\zeta\in
{\mathbb{P}}^1$, we have a~map
\begin{gather*}
\Psi_\zeta: \ {\mathcal{M}}_d\to Z_{\zeta}^{[d]},
\qquad
C\mapsto C_\zeta,
\end{gather*}
to the Douady space of $0$-dimensional complex subspaces of $Z_\zeta$ of length $d$.
This map describes the complex structures of ${\mathcal{M}}_d$:
\begin{Proposition}\label{Hilb}
The map $\Psi_\zeta:({\mathcal{M}}_d,I_\zeta)\to Z_{\zeta}^{[d]} $ is holomorphic and, provided that $\dim Z_\zeta=\dim
Z-1$, an unramified covering of an open subset of the smooth locus of $Z_{\zeta}^{[d]}$.
\end{Proposition}

\begin{proof}
We show f\/irst that the map $\Psi_\zeta:{\mathcal{X}}_d\to Z_{\zeta}^{[d]}$, $\Psi_\zeta(C)= C_{\zeta}$, is holomorphic.
Consider the subspace $T$ of ${\mathcal{X}}_d\times Z_{\zeta}$ def\/ined as
$
T=\{(C,z)\in {\mathcal{X}}_d\times Z_{\zeta};\: z\in C\}$.
Since the projection $\nu:T\to {\mathcal{X}}_d$ is f\/inite
and every f\/ibre has the same length, $\nu$ is f\/lat.
It is also clearly proper.
Hence, the universal property of the Douady space of $Z_{\zeta}$ implies that there is a~unique holomorphic map
$f:{\mathcal{X}}_d\to {\mathcal{D}}(Z_{\zeta})$ such that $T={\mathcal{X}}_d\times_{{\mathcal{D}}(Z_{\zeta})}Y_\zeta$,
where $Y_\zeta$ is the universal family on ${\mathcal{D}}(Z_\zeta)$.
It follows that $f=\Psi_\zeta$, and, hence, $\Psi_\zeta$ is holomorphic on ${\mathcal{X}}_d$.
In addition, this map factors locally through ${\mathcal{X}}_d/Q_\zeta\simeq ({\mathcal{M}}_d,I_\zeta)$, and hence
$\Psi_\zeta$ is holomorphic on~$({\mathcal{M}}_d,I_\zeta)$.

 For $\zeta_1\neq\zeta_2\in {\mathbb{P}}^1$ with $\dim Z_{\zeta_1}=\dim Z_{\zeta_2}=\dim Z-1$, consider now the map
\begin{gather*}
\Phi: \ {\mathcal{X}}_d\ni C\mapsto \bigl( C_{\zeta_1}, C_{\zeta_2}\bigr)\in Z_{\zeta_1}^{[d]}\times Z_{\zeta_2}^{[d]}.
\end{gather*}
An argument analogous to the one given above shows that $\Phi$ is holomorphic.
Since
\[
H^\ast(C,{\mathcal{N}}_{C/Z}(-2))=0,
\]
the dif\/ferential of $\Phi$ is injective at any point of ${\mathcal{X}}_d$.
Thus $\Phi$ is an immersion and, since ${\mathcal{X}}_d$ and $ Z_{\zeta_1}^{[d]}\times Z_{\zeta_2}^{[d]}$ have the same
dimension, $\Phi$ is a~local dif\/feomorphism.
In particular $\Phi$ maps to the smooth locus of~$ Z_{\zeta_1}^{[d]}\times Z_{\zeta_2}^{[d]}$ and is a~covering of its
image.
Taking~$\zeta_2$ to be antipodal to $\zeta_1$, we conclude that~${\mathcal{M}}_d$ is a~covering of an open subset of the
smooth locus of~$ Z_{\zeta_1}^{[d]}$.
\end{proof}

\begin{Remark}\label{twistor}
We can now describe the usual twistor space of ${\mathcal{M}}_d$ as long as $\pi:Z\to {\mathbb{P}}^1$ is f\/lat (e.g.,
submersion of smooth manifolds).
Consider the relative Douady space ${\mathcal{D}}^{[d]}_\pi(Z)$ of f\/inite subspaces of length $d$ in each f\/ibre
$Z_\zeta$~\cite{Pour}, and let $Z_d$ be its open subset consisting of the smooth locus in each f\/ibre.
Since $\pi$ is f\/lat, $Z_d$ is a~manifold.
It has an induced real structure and a~canonical projection $\tilde\pi:Z_d\to {\mathbb{P}}^1$.
Each $C$ in ${\mathcal{M}}_d$ corresponds to a~section $s_C$ of $\tilde\pi:Z_d\to {\mathbb{P}}^1$, $s_C(\zeta)=C_\zeta$,
and the normal bundle ${\mathcal{N}}$ of $s_C({\mathbb{P}}^1)$ in $Z_d$ is isomorphic to $\pi_\ast {\mathcal{N}}_{C/Z}$~-- a~locally free sheaf of rank $2dn$, where $n=\dim Z-1$.
It follows that ${\mathcal{N}}\otimes{\mathcal{O}}_{{\mathbb{P}}^1}(-2)=\pi_\ast \bigl({\mathcal{N}}_{C/Z}(-2)\bigr)$
has no cohomology, and, consequently, ${\mathcal{N}}\simeq \bigoplus{\mathcal{O}}(1)$.
\end{Remark}

\begin{Remark}
The above results can be viewed as follows.
Start with a~hypercomplex mani\-fold~$M$.
Its twistor space is a~smooth mani\-fold~$Z$ equipped with a~projection to~${\mathbb{P}}^1$ and an antiholomorphic
involution which covers the antipodal map on~${\mathbb{P}}^1$.
We obtain, for each~$d$, a~hypercomplex mani\-fold~${\mathcal{M}}_d$, which is biholomorphic, with respect to each complex
structure~$I_\zeta$, to a~discrete covering of an open subset of the smooth locus of the Douady space~$M^{[d]}$ of~$d$
points in~$(M,I_\zeta)$.
Of course, this open subset (and the manifold ${\mathcal{M}}_d$ itself) could be empty.
%\label{X_d}
\end{Remark}

We shall now show that ${\mathcal{M}}_d$ is hyperk\"ahler if $M$ was and $\dim M=4$.
First of all, let us def\/ine symplectic forms in the context we shall need them.
\begin{Definition}
Let $X$ be a~complex space, ${\mathcal{F}}$ a~coherent sheaf, and $L$ a~line a~bundle on~$X$.
An {\em $L$-valued symplectic form} on ${\mathcal{F}}$ is a~homomorphism $\Lambda^2{\mathcal{F}}\to L$ such that the
associated homomorphism ${\mathcal{F}}\to {\mathcal{F}}^\ast\otimes_{{\mathcal{O}}_X} L$ is an isomorphism.
\end{Definition}

Now recall that kernel of the map $d\pi:{\mathcal{T}}_Z\to T{\mathbb{P}}^1$ is called the {\em vertical tangent sheaf}
and is denoted by ${\mathcal{T}}_{Z/{\mathbb{P}}^1}$.
To obtain a~hyperk\"ahler metric we need an ${\mathcal{O}}_Z(2)$-valued symplectic form $\omega$ on
${\mathcal{T}}_{Z/{\mathbb{P}}^1}$.
This form needs to be compatible with the real structure $\sigma$ in the following sense~\cite{HKLR}.
The line bundle ${\mathcal{O}}_Z(1)$ has a~canonical quaternionic structure (i.e.\
an antilinear isomorphism with square $-1$) covering $\sigma$ on $Z$.
Since $\sigma$ induces a~real structure on ${\mathcal{T}}_{Z/{\mathbb{P}}^1}$, we obtain a~quaternionic structure on
${\mathcal{T}}_{Z/{\mathbb{P}}^1}(-1)$.
The form $\omega$ induces a~usual (i.e.\
${\mathcal{O}}$-valued) symplectic form on ${\mathcal{T}}_{Z/{\mathbb{P}}^1}(-1)$, and we say that {\em $\omega$ is
compatible with $\sigma$} if $\omega(\sigma^\ast s,\sigma^\ast t)=\overline{\omega(s,t)}$ for local sections~$s$,~$t$ of~${\mathcal{T}}_{Z/{\mathbb{P}}^1}(-1)$.
We have
\begin{Theorem}
Let $Z$ be a~complex manifold of dimension $3$ equipped with a~holomorphic submersion \mbox{$\pi:Z{\to} {\mathbb{P}}^1$} and an antiholomorphic
involution $\sigma$ covering the antipodal map on ${\mathbb{P}}^1$.
In addition, suppose that we are given a~$\sigma$-compatible ${\mathcal{O}}_Z(2)$-valued symplectic form on the vertical
tangent bundle~${\mathcal{T}}_{Z/{\mathbb{P}}^1}$, which induces a~symplectic structure in the usual sense on each fibre~$Z_\zeta$.
Then the hypercomplex manifold ${\mathcal{M}}_d$ defined in Theorem~{\rm \ref{hc}} has a~canonical pseudo-hyperk\"ahler
metric.
\end{Theorem}

\begin{Remark}
The signature of this metric can vary between dif\/ferent connected components of ${\mathcal{M}}_d$.
\end{Remark}
\begin{proof}
The arguments of Beauville~\cite{Beau} show that if $M$ is a~complex surface with a~holomorphic symplectic form, then
 $M^{[d]}$ has an induced holomorphic symplectic form.
Applying this construction f\/ibrewise to $Z$, we obtain an ${\mathcal{O}}(2)$-valued f\/ibrewise symplectic form on the
twistor space $Z_d$ of ${\mathcal{M}}_d$ described in Remark~\ref{twistor}.
This form is still compatible with the induced real structure and hence, owing to~\cite[Theorem 3.3]{HKLR}, it gives
a~pseudo-hyperk\"ahler metric on ${\mathcal{M}}_d$.
\end{proof}

\begin{Remark}
Nash~\cite{Nash} gives a~dif\/ferent construction of hyperk\"ahler metrics on moduli spaces of ${\rm SU}(2)$-monopoles, which
works as long as  ${\mathcal{M}}_d$ is replaced by its open subset of l.c.i.\
curves (so that the normal sheaf is locally free and the Serre duality can be applied to~it).
\end{Remark}

\begin{Remark}
Similarly, if $Z$ has an ${\mathcal{O}}(n)$-valued volume form on ${\mathcal{T}}_{Z/{\mathbb{P}}^1}$, then so does
$Z_d$, and the holonomy group of ${\mathcal{M}}_d$ reduces to ${\rm SL}(dn,{\mathbb{H}})$.
\end{Remark}

\section{Examples}

\subsection{Monopoles and the generalised Legendre transform}

In the case of ${\rm SU}(2)$-monopoles, one starts with the twistor space $Z$ of $S^1\times {\mathbb{R}}^3$, i.e.\
the total space of certain line bundle $L^2$ over $T{\mathbb{P}}^1$ without the zero section.
Write $p:Z\to T{\mathbb{P}}^1$ for the projection.
A curve $C$ in $Z$ corresponds to a~curve $p(C)$ in $T{\mathbb{P}}^1$ such that $L^2|_{p(C)}$ is trivial.
This is the condition satisf\/ied by spectral curves of magnetic monopoles, and Nash~\cite{Nash} shows that if $S$ is
a~spectral curve of a~monopole, then its lift $C$ to $L^2$ satisf\/ies additionally $H^\ast({\mathcal{N}}_{C/Z}(-2))=0$.
Thus, the moduli space of monopoles of charge $d$ (i.e.\
those for which $S$ is a~curve of degree $d$) is a~connected component of the manifold ${\mathcal{M}}_d$ def\/ined in the
previous section.
Already in this case one has to include singular curves (although not nonreduced ones).

More general hyperk\"ahler metrics were considered in~\cite{Bie}, as examples of the generalised Legendre construction
of Lindstr\"om and Ro\v{c}ek~\cite{LR}.
Many of these can be put in this framework, i.e.~$Z$ is the total space of a~line bundle over a~complex surface, or, more generally, the total space of a~holomorphic
principal bundle over a~complex manifold f\/ibring over~${\mathbb{P}}^1$, or the projectivisation of a~vector bundle over
a~complex manifold f\/ibring over~${\mathbb{P}}^1$.
The last situation is relevant, for example, for ALF gravitational instantons of type~$D_k$~\cite{CH}.

According to Proposition~\ref{Hilb}, the complex structures of such hypercomplex manifolds are always those of (covering
of) open subsets of the Douady spaces of $0$-dimensional subspaces of f\/ibres of $Z$.
In the case when $Z$ is the total space of a~${\mathbb C}^\ast$-bundle over a~complex surface $\Sigma\to
{\mathbb{P}}^1$, one can describe these open subsets more precisely, as in~\cite[Chapter~6]{AH}.
Let $p:Z\to \Sigma$ be the projection.
Then $({\mathcal{M}}_d, I_\zeta)$ is biholomorphic to the open subset of $(Z_\zeta)^{[d]}_p$, where the subscript~$p$
denotes $0$-dimensional subspaces $D$ such that $p_\ast({\mathcal{O}}_D)$ is a~cyclic
${\mathcal{O}}_{\Sigma_\zeta}$-sheaf.

\subsection{Projective curves}

The twistor space $Z$ of the f\/lat ${\mathbb{R}}^4$ is the total space of ${\mathcal{O}}(1)\oplus {\mathcal{O}}(1)$ on
${\mathbb{P}}^1$.
We can view it as ${\mathbb{P}}^3-{\mathbb{P}}^1$, so that curves in $Z$ are curves in the projective space not
intersecting a~f\/ixed projective line.
In addition, the real structure of $Z$ extends to the real structure $\sigma$ of ${\mathbb{P}}^3$ (which is the twistor
space of $S^4$), so we look for $\sigma$-invariant curves in ${\mathbb{P}}^3$.

Consider f\/irst a~rational curve, i.e.\
an embedding $f:{\mathbb{P}}^1\to {\mathbb{P}}^3$ given by homogeneous polynomials of degree~$d$.
For $d=2$, such a~curve is an intersection of a~line and a~quadric, so its normal bundle is isomorphic to
${\mathcal{O}}_{{\mathbb{P}}^1}(2)\oplus {\mathcal{O}}_{{\mathbb{P}}^1}(4)$ and it does not satisfy the condition
$H^\ast({\mathcal{N}}(-2))=0$.
For $d\geq 3$, however, the normal bundle of a~generic rational curve of deg\-ree~$d$ splits as~$
{\mathcal{O}}_{{\mathbb{P}}^1}(2d-1)\oplus {\mathcal{O}}_{{\mathbb{P}}^1}(2d-1)$~\cite{EV,S?}, while the restriction of
${\mathcal{O}}_Z(1)= {\mathcal{O}}_{{\mathbb{P}}^3}(1)|_Z$ is isomorphic to $ {\mathcal{O}}_{{\mathbb{P}}^1}(d)$.
Thus a~generic rational curve of degree $d\geq 3$ satisf\/ies the condition $H^\ast({\mathcal{N}}_{C/Z}(-2))=0$ and,
consequently, the parameter space of such curves, which are $\sigma$-invariant and avoid a~f\/ixed line, is
a~$4d$-dimensional (pseudo)-hyperk\"ahler manifold.
We shall see shortly, that for $d=3$, i.e.\
for twisted normal cubics, the resulting metric is f\/lat.
This, however, is a~rather special case and we do not know what to expect in the general case.
We observe that the action of ${\rm SO}(3)$ on ${\mathbb{P}}^1$ induces an isometric action rotating the complex structures,
and so all complex structures are equivalent.
In fact, we expect that these hyperk\"ahler manifolds are cones over $3$-Sasakian manifolds.

For higher genera, it is known~\cite{EH} that the parameter space $H_{d,g}$ of space curves with degree~$d$ and genus~$g$ contains smooth curves with $H^\ast({\mathcal{N}}_{C/{\mathbb{P}}^3}(-2))=0$ for any~$d$ greater than or equal to
some $D(g)$ (e.g., $D(0)=3$ and $D(1)=5$).
As soon as $H_{d,g}$ contains also a~{\em $\sigma$-invariant} smooth curve with
$H^\ast({\mathcal{N}}_{C/{\mathbb{P}}^3}(-2))=0$, we obtain a~natural pseudo-hyperk\"ahler structure on a~submanifold of
$H_{d,g}$.

We shall now show that the resulting metric on $H_{3,0}$ (and more generally on moduli spaces of ACM (arithmetically
Cohen--Macaulay) curves admitting a~linear resolution) is f\/lat.
Although the metric itself is not interesting, it is still an instructive example, which shows, in particular, that if
we want to have a~shot at completeness of the metric, we cannot avoid including very singular and nonreduced complex
subspaces of $Z$.

\subsection{ACM curves with a~linear resolution}

We consider curves $C$, the structure sheaf of which admits a~free resolution of the form
\begin{gather}
0\to {\mathcal{O}}_{{\mathbb{P}}^3}(-r-1)^r\stackrel{\phi_2}{\longrightarrow} {\mathcal{O}}_{{\mathbb{P}}^3}(-r)^{r+1}
\stackrel{\phi_1}{\longrightarrow} {\mathcal{O}}_{{\mathbb{P}}^3} \longrightarrow {\mathcal{O}}_C\to 0.
\label{res}
\end{gather}
This means that ${\mathcal{I}}_C$ is def\/ined by simultaneous vanishing of the $r\times r$ minors of the linear matrix
$\phi_2$.
If $C$ is smooth, then its degree is $d=\frac{1}{2}r(r+1)$ and its genus is equal to $g=\frac{1}{6}(r-1)(r-2)(2r+3)$.

A complex subspace $C$ with a~resolution~\eqref{res} is automatically equidimensional and Cohen--Macaulay~\cite[Theorem~18.18]{Eis1}, and the Douady space (i.e.\
Hilbert scheme) is smooth at $C$~\cite[Co\-rol\-lary~8.10]{Hartsh}.
Furthermore, Ellia~\cite{Ellia} has shown that every such $C$ satisf\/ies \mbox{$H^\ast(C,{\mathcal{N}}_{C/Z}(-2))=0$}.
For completeness (and to remove the unnecessary assumption of projective normality) let us reproduce his argument.
\begin{Lemma}[\protect{\cite{Ellia}}]
Let $C$ be a~subscheme of ${\mathbb{P}}^3$ with a~resolution~\eqref{res}.
Then $H^\ast({\mathcal{N}}_{C/Z}(-2))=0$.
\end{Lemma}

\begin{proof}
We can rewrite~\eqref{res} as a~resolution of the ideal sheaf of $C$:
\begin{gather}
0\to {\mathcal{O}}_{{\mathbb{P}}^3}(-r-1)^r\stackrel{\phi_2}{\longrightarrow} {\mathcal{O}}_{{\mathbb{P}}^3}(-r)^{r+1}
\stackrel{\phi_1}{\longrightarrow} {\mathcal{I}}_C\to 0.
\label{res2}
\end{gather}
We have $H^i({\mathcal{N}}_{C/Z}(-2))=\operatorname{Ext}^{i+1}({\mathcal{I}}_C,{\mathcal{I}}_C(-2))$, $i=0,1$.
Applying $\operatorname{Hom}(-,{\mathcal{I}}_C(-2))$ to~\eqref{res2} and using the isomorphism
$\operatorname{Ext}^i({\mathcal{O}}_{{\mathbb{P}}^3}(k), {\mathcal{I}}_C(-2))\simeq H^i({\mathcal{I}}_C(-2-k))$, $k\in
{\mathbb{Z}}$, we obtain the exact sequence
\begin{gather*}
\begin{split}
& H^0({\mathcal{I}}_C(r-1))^r\to \operatorname{Ext}^1({\mathcal{I}}_C,{\mathcal{I}}_C(-2))\to
H^1({\mathcal{I}}_C(r-2))^{r+1}
\\
& \hphantom{H^0({\mathcal{I}}_C(r-1))^r}{}
\to H^1({\mathcal{I}}_C(r-1))^r\to \operatorname{Ext}^2({\mathcal{I}}_C,{\mathcal{I}}_C(-2))\to 0.
\end{split}
\end{gather*}
On the other hand, tensoring the resolution~\eqref{res2} with $ {\mathcal{O}}_{{\mathbb{P}}^3}(r-1)$ and with $
{\mathcal{O}}_{{\mathbb{P}}^3}(r-2)$ shows that ${\mathcal{I}}_C(r-1)$ and ${\mathcal{I}}_C(r-2)$ are acyclic.
\end{proof}

Let us now choose a~$\sigma$-invariant $2$-dimensional linear system $H$ in ${\mathcal{O}}_{{\mathbb{P}}^3}$ with base
$B$, and let $Z={\mathbb{P}}^3-B$ be the twistor space of ${\mathbb{R}}^4$.
The projection $\pi:Z\to {\mathbb{P}}^1$ is def\/ined by $H$, and it is automatically f\/inite-to-one on any projective
subscheme which is contained in ${\mathbb{P}}^3-B$, since a~$1$-dimensional intersection with $\pi^{-1}(\zeta)$ will
also intersect $B$.

Thus, the Douady space of those $\sigma$-invariant $C$, which admit a~resolution of the form~\eqref{res} and do not
intersect a~f\/ixed ${\mathbb{P}}^1$, is a~(pseudo)-hyperk\"ahler manifold $X_r$ of dimension $2r(r+1)$.
We shall now show that $X_r$ is the f\/lat ${\mathbb{H}}^{r(r+1)/2}$.
Moreover, the natural biholomorphism of Proposition~\ref{Hilb} identif\/ies~$X_r$ with an open subset of the Hilbert
scheme $({\mathbb C}^2)^{[r(r+1)/2]}$ consisting of $0$-dimensional subspaces of length~$r$ which are not subschemes of
any plane curve of degree~${<}r$.

Let $x_1,\dots,x_4$ be homogeneous coordinates on ${\mathbb{P}}^3$ and $\pi:[x_1,\dots,x_4]\mapsto [x_3,x_4]$ the chosen
projection onto ${\mathbb{P}}^1$.
Thus $B=\{[x_1,x_2,0,0]\}$ is the base of the linear system and ${\mathbb{P}}^3-B$ is our twistor space $Z$.
The real structure $\sigma$ is given by
\begin{gather*}
\sigma: \ [x_1,x_2,x_3,x_4]\mapsto [-\overline{x}_2,\overline{x}_1,-\overline{x}_4,\overline{x}_3].
\end{gather*}
Let $C$ be a~curve def\/ined by~\eqref{res}, and write, with respect to some bases,
\begin{gather*}
\phi_2(x_1,x_2,x_3,x_4)=\sum\limits_{i=1}^4 A_ix_i,
\qquad
A_i\in \operatorname{Mat}_{r+1,r}({\mathbb C}).
\end{gather*}
Such a~$C$ does not intersect $B$ if $\phi_2$ restricted to
$x_3=x_4=0$ has rank $r$ for all non-zero $(x_1,x_2)$.
The involution $\sigma$ induces a~quaternionic structure $\sigma$ on linear forms and such a~curve~$C$ is
$\sigma$-invariant as soon as
\begin{itemize}\itemsep=0pt
\item
for $r$ even, the set of columns of $\phi_2$ is $\sigma$-invariant;
\item
for $r$ odd, the set of rows of $\phi_2$ is $\sigma$-invariant.
\end{itemize}

We shall now describe the intersection $C_{\zeta}$ of $C$ with a~f\/ibre $\pi^{-1}(\zeta)$.
The map $\phi_2$ restricted to the projective plane~$\overline{Z_\zeta}$ is still injective, otherwise~$C_\zeta$ is
$2$-dimensional and intersects~$B$.
Thus~$C_\zeta$ also has a~free resolution of the form~\eqref{res}:
\begin{gather}
0\to {\mathcal{O}}_{{\mathbb{P}}^2}(-r-1)^r\stackrel{\psi_2}{\longrightarrow} {\mathcal{O}}_{{\mathbb{P}}^2}(-r)^{r+1}
\stackrel{\psi_1}{\longrightarrow} {\mathcal{O}}_{{\mathbb{P}}^2} \longrightarrow {\mathcal{O}}_{C_\zeta}\to 0.
\label{res3}
\end{gather}
Let us now write $\zeta=[1,t]$, so that $\psi_2(x_1,x_2,x_3)=A_1x_1+A_2x_2+(A_3+tA_4)x_3$.
Recall that the condition that $C$ does not intersect $B$ is equivalent to $A_1x_1+A_2x_2$ having rank $r$ for any
$[x_1,x_2]\in {\mathbb{P}}^1$.

\begin{Lemma}
If $A_1$, $A_2$ are two $(r+1)\times r$ complex matrices such that $ A_1x_1+A_2x_2$ is injective for every $(x_1,x_2)\in
{\mathbb C}^2\backslash\{0\}$, then $ A_1x_1+A_2x_2$ belongs to the $($open$)$ ${\rm GL}(r+1,{\mathbb C})\times {\rm GL}(r,{\mathbb
C})$-orbit of $ Sx_1 +Tx_2$ where
\begin{gather*}
S_{ij}=
\begin{cases}
1 & \text{if}\quad i=j,
\\
0 & \text{if}\quad i\neq j,
\end{cases}
\qquad
T_{ij}=
\begin{cases}
1 & \text{if}\quad i=j+1,
\\
0 & \text{if}\quad i\neq j+1.
\end{cases}
\end{gather*}
Moreover, the stabiliser of $(A_1,A_2)$ in ${\rm GL}(r+1,{\mathbb C})\times {\rm GL}(r,{\mathbb C})$ is the central subgroup
\begin{gather*}
\Delta=\big\{(z{\rm Id},z^{-1}{\rm Id}); z\in {\mathbb C}^\ast\big\}.
\end{gather*}
\end{Lemma}

\begin{proof}
We can appeal to Kronecker's theory of minimal indices~\cite[Chapter~XII]{Grant}.
The assumption implies that the pencil $A_1+\lambda A_2$ has no minimal indices for columns and no elementary divisors.
Thus it can have only minimal indices for rows, and therefore it lies in the ${\rm GL}(r+1,{\mathbb C})\times {\rm GL}(r,{\mathbb
C})$-orbit of a~block quasi-diagonal matrix built out of blocks of the form as in the statement of the lemma.
Such a~matrix cannot, however, have size $(r+1)\times r$, unless there is only one block.
The statement about the stabiliser is then straightforward.
\end{proof}

We can now use the action of ${\rm GL}(r+1,{\mathbb C})\times {\rm GL}(r,{\mathbb C})$ in order to make $Sx_1+Tx_2$
$\sigma$-invariant (and f\/ixed).
It is then easy to see that, given an $(r+1)\times r$ matrix $\tilde A_3$, we can f\/ind a~unique~$A_3$,~$A_4$ so that
$\tilde A_3=A_3+tA_4$ and $A_3x_3+A_4x_4$ is $\sigma$-invariant.
It follows that the map which sends~$C_\zeta$ represented by $(S,T,\tilde A_3)$ to $\tilde A_3$ is an isomorphism
between the twistor space $Z_r$ of $X_r$ (cf.\
Remark~\ref{twistor}) and the total space of ${\mathbb C}^{r(r+1)}\otimes {\mathcal{O}}(1)$, i.e.\
the twistor space of the f\/lat ${\mathbb{H}}^{r(r+1)/2}$ (possibly, although unlikely, with a~non-Euclidean signature).

It remains to identify the complex structure of $X_r$ as an open subset of the Hilbert scheme $({\mathbb
C}^2)^{[r(r+1)/2]}$.
Recall that the Hilbert scheme $({\mathbb C}^2)^{[d]}$ of $d$ points in ${\mathbb C}^2$ has a~natural stratif\/ication by
the Hilbert function $H:{\mathbb{N}}\to {\mathbb{N}}$, with~$H(k)$ equal to $(k+1)(k+2)/2$ minus the dimension of the
vector space of plane curves of degree~$k$ containing~$X$~\cite{Eis1}.
The Hilbert function can be computed from any free resolution, and in our case we obtain
\begin{gather*}
H(k)=
\begin{cases}
\dfrac{(k+1)(k+2)}{2} &\text{if}\quad k<r,
\vspace{1mm}\\
\dfrac{r(r+1)}{2} & \text{if}\quad k\geq r.
\end{cases}
%\label{H(k)}
\end{gather*}
This means that $C_\zeta$ does not lie on any plane curve of degree~${<}r$.
Conversely, if a~$D\in ({\mathbb{P}}^2)^{[d]}$ does not lie on any curve of degree~${<}r$, then its ideal must be
generated by forms of degree~${\geq}r$.
Comparing dimensions and using the Hilbert--Burch theorem, as in~\cite[\S~20.4.1]{Eis1}, shows that the minimal
resolution of $D$ is of the form~\eqref{res3}.
In other words, the natural biholomorphism of Proposition~\ref{Hilb} identif\/ies~$X_r$ with the open stratum of
$({\mathbb C}^2)^{[r(r+1)/2]}$, consisting of $0$-dimensional subspaces which do not lie on any plane curve of degree
smaller than~$r$.

\section{Induced vector bundles and their tangent spaces}%\label{vb}

We aim to def\/ine canonical connections on vector bundles over a~hypercomplex manifold, obtained from higher degree
curves in $Z$, directly in terms of objects def\/ined on~$Z$, i.e.\
without passing to the usual twistor space of a~hypercomplex manifold.
In this section, we shall consider a~general Douady space~${\mathcal{D}}(Z)$ (not necessarily a~hypercomplex manifold)
and describe total spaces of vector bundles arising on~${\mathcal{D}}(Z)$ as Douady spaces of some other, canonically
def\/ined, complex spaces.

Let $Z$ be a~complex space and ${\mathcal{D}}(Z)$ its Douady space.
The double f\/ibration
\begin{gather*}
{\mathcal{D}}(Z)\stackrel{\nu}{\longleftarrow} Y \stackrel{\mu}{\longrightarrow} Z
%\label{seq0}
\end{gather*}
allows one to transfer holomorphic data from $Z$ to ${\mathcal{D}}(Z)$ or to its subsets.
Let $M$ be an open connected subset of ${\mathcal{D}}(Z)$, such that every $X\in M$ is regularly embedded (cf.\
Section~\ref{back}) and satisf\/ies $h^1({\mathcal{N}}_{X/Z})=0$.
In particular, $M$ is a~manifold and ${\mathcal{N}}_{X/Z}$ a~locally free sheaf.
We are interested in vector bundles on $M$ obtained from vector bundles on $Z$.
Let ${\mathcal{E}}$ be a~vector bundle on $Z$.
The sheaf $\widehat{\mathcal{E}}=\nu_\ast\mu^\ast {\mathcal{E}}$ is locally free on~$M$ if the function $X\mapsto
h^0({\mathcal{E}}|_{X})$ is constant on~$M$.
In this case $\widehat{\mathcal{E}}$ is called the {\em vector bundle induced by ${\mathcal{E}}$}~\cite{HM}.

Let $\pi:E\to M$ be a~vector bundle on a~manifold $M$.
One way of def\/ining connections on $E$ is to split the canonical exact sequence
\begin{gather}
0\to V\to TE\to \pi^\ast(TM)\to 0
\label{seq1.5}
\end{gather}
of vector bundles on $E$, where $V$ denotes the vertical bundle $\ker d {\pi}$.
As the f\/irst step in def\/ining canonical connections, we want to show that for bundles $\widehat{\mathcal{E}}$ (with
${\mathcal{E}}$ satisfying a~cohomological condition) the above sequence is induced from objects on $Z$.

Let ${\mathcal{E}}$ be a~vector bundle on $Z$ such that $h^0({\mathcal{E}}|_X)$ is constant on~$M$ and
$h^1({\mathcal{E}}|_{X})=0$ for every~$X$ in~$M$.
A point of the total space of $E=\widehat{{\mathcal{E}}}$ corresponds to a~pair, consisting of a~regularly embedded
compact subspace $X$ in $Z$ and a~section $s$ of ${\mathcal{E}}|_{X}$.
The image space $s(X)$ is well def\/ined and is a~regularly embedded compact subspace of the total space\footnote{The
reason for the inconsistence in writing $|{\mathcal{E}}|$ for the total space of ${\mathcal{E}}$ is that
${\mathcal{D}}({\mathcal{E}})$ has a~dif\/ferent meaning: it parameterises coherent quotients of ${\mathcal{E}}$ with
compact support.} $|{\mathcal{E}}|$ of ${\mathcal{E}}$.
The normal sheaf ${\mathcal{N}}_{s(X)}$ of $s(X)$ in $|{\mathcal{E}}|$ f\/its into the exact sequence
\begin{gather}
0\to \pi^\ast {\mathcal{E}}|_{X}\to {\mathcal{N}}_{s(X)}\to \pi^\ast {\mathcal{N}}_{X/Z}\to 0,
\label{seq1}
\end{gather}
of sheaves on $s(X)$.
Since $h^1({\mathcal{E}}|_{X})=0$, it follows that $h^1({\mathcal{N}}_{s(X)})=0$ and the Douady space~${\mathcal{D}}\bigl(|{\mathcal{E}}|\bigr)$ is smooth at $s(X)$ of dimension $h^0( {\mathcal{N}}_{s(X)})=\dim M
+\operatorname{rank} E$.
Since pairs $(X^\prime,s^\prime)$ def\/ine a~submanifold of ${\mathcal{D}}\bigl(|{\mathcal{E}}|\bigr)$ also of dimension
$\dim M +\operatorname{rank} E$, it follows that $E$ is an open subset of the smooth locus of
${\mathcal{D}}\bigl(|{\mathcal{E}}|\bigr)$.
Therefore the tangent space $T_eE$ at $e=(X,s)$ is cano\-nical\-ly identif\/ied with $H^0({\mathcal{N}}_{s(X)})$ and the
sequence~\eqref{seq1} induces the canonical sequence~\eqref{seq1.5}.
We can phrase these considerations more precisely, if less transparently, as follows.

\begin{Theorem}
Let $Z$ be a~complex space and~$M$ an open connected subset of the smooth locus of~${\mathcal{D}}(Z)$, such that every
$X\in M$ is regularly embedded.
Let ${\mathcal{E}}$ be a~vector bundle on~$Z$ such that $h^0({\mathcal{E}}|_X)$ is constant on~$M$ and
$h^1({\mathcal{E}}|_{X})=0$ for every~$X$ in~$M$, and let $E=\widehat{\mathcal{E}}$ be the induced vector bundle on~$M$.
Then~$E$ is identified with an open subset of the smooth locus of~${\mathcal{D}}\bigl(|{\mathcal{E}}|\bigr)$, again
consisting of regularly embedded subspaces.
Furthermore, if~$\tilde Y$ is the universal family on~${\mathcal{D}}\bigl(|{\mathcal{E}}|\bigr)$ and~$p:\tilde Y\to Y$,
$\tilde \nu:\tilde Y\to E$ are canonical projections, then the normal sheaf ${\mathcal{N}}_{\tilde Y}$ of $\tilde Y$ in
$E\times |{\mathcal{E}}|$ fits into the commutative diagram:
\begin{gather}
\begin{CD}
0 @>>> \tilde\nu_\ast p^\ast\mu^\ast {\mathcal{E}} @>>> \tilde\nu_\ast{\mathcal{N}}_{\tilde Y} @>>> \tilde\nu_\ast
p^\ast {\mathcal{N}}_Y @>>> 0
\\
@.
@V V\wr V @V V\wr V @V V\wr V @.
\\
0 @>>> V @>>> TE @>>> \pi^\ast(TM) @>>>0
\end{CD}
\end{gather}
\end{Theorem}

\begin{Remark}
The above description of the total space of an induced vector bundle as a~Douady space has the following consequence.
Let $Z$ be as in Theorem~\ref{hc} and let $M$ be an open subset of ${\mathcal{M}}_d$ consisting of regularly embedded
curves.
Let ${\mathcal{E}}$ be a~vector $\sigma$-bundle on $Z$ (i.e.\
there is an antilinear involution on ${\mathcal{E}}$ covering $\sigma$ on $Z$) such that $h^0({\mathcal{E}}|_C)$ is
constant on $M$ and ${\mathcal{E}}(-2)$ is acyclic on any $C\in M$.
Then the $\sigma$-invariant part of the total space of $\widehat {\mathcal{E}}$ has a~natural hypercomplex structure.
Indeed, it is enough to tensor~\eqref{seq1} by ${\mathcal{O}}(-2)$ and apply the results of Section~\ref{hgc}.
This result is in the same spirit as~\cite[Theorem~7.2]{Sal}.
\end{Remark}

\section{Hyperholomorphic connections}

We return to hypercomplex manifolds and consider again the situation from Section~\ref{hgc}.
As in the last section, we restrict ourselves to the subset of regularly embedded curves, for which the construction is
more transparent.
Thus~$Z$ is equipped with a~holomorphic surjection onto~${\mathbb{P}}^1$ and~$M$ is an open connected subset of~${\mathcal{D}}(Z)$ consisting of regularly embedded compact subspaces~$C$, such that~$\pi|_C$ is f\/inite-to-one and
$H^\ast (C,{\mathcal{N}}_{C/Z}(-2))=0$ (cf.\
Lemma~\ref{ooo}).
Recall that~$M$ is then a~complex version of a~hypercomplex manifold, i.e.\
$TM\simeq E\otimes {\mathbb{C}}^2$, and for every $\zeta\in {\mathbb{P}}^1$ the distribution $Q_\zeta=E\otimes l$, where
$[l]=\zeta$, is integrable.
We now use the results of the previous section to def\/ine canonical connections on certain induced vector bundles on~$M$.

\begin{Proposition}
Let ${\mathcal{E}}$ be a~vector bundle on $Z$ such that $C\to h^0({\mathcal{E}}|_C)$ is constant on $M$ and
${\mathcal{E}}(-1)|_{C}$ is acyclic for every~$C\in M$.
Then the induced bundle $\widehat{{\mathcal{E}}}$ on $M$ is equipped with a~canonical linear connection $\nabla$, which
has the following property: for any $\zeta\in {\mathbb{P}}^1$ and any $m\in M$, if $u$ is a~local section of
$\widehat{{\mathcal{E}}}$ which satisfies $du(X)=0$ for every $X\in Q_\zeta|_m$, then $\nabla_X u=0$ for every such $X$.

In particular, if $Z$ is also equipped with an antiholomorphic involution $\sigma$ covering the antipodal map, so that
$M^\sigma$ is a~hypercomplex manifold, then $\nabla$ on $\widehat{{\mathcal{E}}}|_{M^\sigma}$ is hyperholomorphic, i.e.\
$\nabla^{0,1}=\bar{\partial}$ for every complex structure $I_\zeta$.
%\label{induced}
\end{Proposition}

\begin{Remark}
In the case when $Z$ is the usual twistor space of a~hypercomplex manifold, i.e.\
the curves $C$ have genus $0$, then the condition that ${\mathcal{E}}(-1)|_{C}$ is acyclic is equivalent to
${\mathcal{E}}|_C$ to being trivial, and we recover the well-known results of~\cite{W} and~\cite{HM}.
In that case, already the sequence~\eqref{seq1} splits.
\end{Remark}

\begin{proof}
Let us write $E$ for $\widehat{{\mathcal{E}}}$ (in the course of the proof we are not going to use the bundle~$E$ from
the decomposition $TM\simeq E\otimes H$), and consider the canonical exact sequence~\eqref{seq1.5}.
An Ehresmann connection is a~splitting of this sequence.
To def\/ine such a~splitting, tensor~\eqref{seq1} by~${\mathcal{O}}(-1)$.
Acyclity of~${\mathcal{E}}(-1)|_{C}$ implies that there is a~canonical isomorphism
\begin{gather*}
\phi: \ H^0({\mathcal{N}}_{C/Z}(-1))\to H^0( {\mathcal{N}}_{s(C)}(-1)).
\end{gather*}
We obtain a~canonical map (where $H=\pi^\ast H^0({\mathbb{P}}^1,{\mathcal{O}}(1))$)
\begin{gather*}
H^0({\mathcal{N}}_{C/Z})\simeq H^0({\mathcal{N}}_{C/Z}(-1))\otimes H\stackrel{\phi\otimes{\rm id}}{\longrightarrow} H^0(
{\mathcal{N}}_{s(C)}(-1))\otimes H \to H^0( {\mathcal{N}}_{s(C)}),
\end{gather*}
which splits~\eqref{seq1.5}.
Since, for any $s\in H$ with $s(\zeta)=0$, the image in ${\mathcal{N}}_{s(C)}$ of $H^0( {\mathcal{N}}_{s(C)}(-1))\otimes
s$ consists of sections which vanish on $s(C)_\zeta$ (here, as in Section~\ref{hgc}, the subscript $\zeta$ denotes the f\/ibre
of the projection to ${\mathbb{P}}^1$), we have a~description of the horizontal subspace at $(C,s)$ as the subspace of
$H^0( {\mathcal{N}}_{s(C)})$ linearly generated by sections vanishing on some $s(C)_\zeta$, $\zeta\in {\mathbb{P}}^1$.
It follows easily that the dif\/ferential of the scalar multiplication $E\to E$ and the dif\/ferential of the addition map
$E\oplus E\to E$ preserve the horizontal subbundles, and, therefore, our Ehresmann connection is linear.

Let us now prove the stated property of this connection.
Let $\gamma(t)$ be a~curve in $M$ tangent to the distribution $Q_\zeta$.
The curves $C_t$ in $Z$ corresponding to $\gamma(t)$ have a~f\/ixed intersection with the f\/ibre $Z_\zeta$.
A horizontal lift $\tilde\gamma(t)$ of $\gamma$ is given by $(C_t,s_t)$ such that $s_t(C_t)\cap {\mathcal{E}}_\zeta$ is
a~f\/ixed complex subspace of the f\/ibre ${\mathcal{E}}_\zeta$.
Thus, for a~section $u$ of $E$, the parallel transport $\tau^t(u(t))= (\tilde C_t,\tilde s_t)$ of $u(\gamma(t))$ along
$\gamma$ to $\gamma(0)$ satisf\/ies
\begin{gather}
\tilde s_t(\tilde C_t)\cap {\mathcal{E}}_\zeta= s_t(C_t)\cap {\mathcal{E}}_\zeta.
\label{parallel}
\end{gather}
Suppose now that $u$ is a~section of $E$ such that $du(X)=0$ for every $X\in Q_\zeta|_m$.
This implies that if $\gamma(t)$ is a~curve in $M$ with $\dot\gamma(0)=X\in Q_\zeta|_m$ and $(\Gamma_t,\psi_t)$ is the
pair {\em curve $+$ section} corresponding to $u(\gamma(t))$, then up to order $1$:{\samepage
\begin{gather*}
\psi_t(\Gamma_t)\cap {\mathcal{E}}_\zeta=\psi_0(\Gamma_0)\cap{\mathcal{E}}_\zeta.
\end{gather*}
Comparing this and~\eqref{parallel}, we conclude that $\nabla_X u=0$ for any $X\in Q_\zeta|_m$.}

The remaining statement of the Proposition (concerning $M^\sigma$) is automatic, since $Q_\zeta$ consists of vectors of
type $(0,1)$ for $I_\zeta$.
\end{proof}

\begin{Remark}
A similar result remains valid for the Grassmann structures def\/ined in Remark~\ref{Grass}.
The above argument produces a~canonical linear connection on a~vector bundle induced from a~vector bundle
${\mathcal{E}}$ on $Z$ with constant $h^0({\mathcal{E}}|_C)$ and vanishing $H^\ast({\mathcal{E}}(-1)|_C)$.
Once again, this connection has the property described in the statement of the above theorem.
\end{Remark}

\begin{Remark}
The method of this and the last section reproduces also another result of Huggett and Merkulov~\cite{HM}.
Let $Z$ be an arbitrary complex space, and $M\subset {\mathcal{D}}(Z)$ as in the previous section.
Let ${\mathcal{E}}$ be a~vector bundle on $Z$ such that $h^0({\mathcal{E}}|_X)$ is constant on $M$ and
\mbox{$h^i({\mathcal{E}}\otimes_X {\mathcal{N}}_{X/Z}^\ast)=0$} for $i=0,1$ and every $X\in M$.
Then there is a~canonical splitting of the sequence~\eqref{seq1}, which produces a~canonical linear Ehresmann connection
as above.
\end{Remark}

\pdfbookmark[1]{References}{ref}
\LastPageEnding

\end{document}